\documentclass{amsart}

\usepackage{amsmath, amsthm, amssymb, graphicx, enumerate, everypage, datetime, getfiledate, lmodern}

\newcounter{forwardtheorem}

\newcounter{citedtheorems}

\newtheorem{defn}{Definition}[section]
\newtheorem{theorem}[defn]{Theorem}
\newtheorem*{theorem-x}{Theorem}
\newtheorem*{theorem-m}{Main Theorem}
\newtheorem*{theorem-n}{Main Theorem}
\newtheorem*{cor-x}{Corollary}
\newtheorem*{lemma-x}{Lemma}
\newtheorem*{concl-x}{Conclusion}
\newtheorem*{claim-x}{Claim}
\newtheorem*{thm-r}{Theorem \ref{concl:sop2-max}}
\newtheorem*{thm-q}{Theorem \ref{theorem:p-t}}
\newtheorem*{claim-s}{Claim \ref{m1-sat}}

\newtheorem{thm-lit}[citedtheorems]{Theorem}
\newtheorem{defn-lit}[citedtheorems]{Definition}

\newtheorem{conv-r}[defn]{Conventions and Remarks}

\newcommand{\xt}{\mathfrak{t}}

\newcommand{\xp}{\mathfrak{p}}

\newcommand{\vp}{\varphi}

\newcommand{\ts}{\mathbf{S}}

\newcommand{\ml}{\mathcal{L}}

\newcommand{\mcp}{\mathcal{P}}

\newcommand{\lost}{\L os' }

\newcommand{\de}{\mathcal{D}}

\newcommand{\blank}{ \textcolor{white}{space}  }

\newcommand{\mcf}{\mathcal{F}}

\newcommand{\br}{\vspace{2mm}}
\newcommand{\sbr}{\vspace{1mm}}

\title{Model theory and ultraproducts}

\author{M. Malliaris}

\address{Department of Mathematics, University of Chicago, 5734 S. University Avenue, Chicago, IL 60637, USA}
\begin{document}

\begin{abstract}
The article motivates recent work on saturation of ultrapowers from a general mathematical point of view. 
\end{abstract}

\maketitle

\br

\subsection*{Introduction} \blank 

\br

In the history of mathematics the idea of the limit has had a remarkable effect on organizing and making visible a certain kind of 
structure. Its effects can be seen from calculus to extremal graph theory. 

At the end of the nineteenth century, Cantor introduced infinite cardinals, 
which  allow for a stratification of size in a potentially much deeper way. 
It is often thought that these ideas, however beautiful, are studied in modern mathematics primarily by set theorists, and that aside from 
occasional independence results, the 
hierarchy of infinities has not profoundly influenced, say, algebra, geometry, analysis or topology, and perhaps is not suited to do so. 

What this conventional wisdom misses is a powerful kind of reflection or transmutation 
arising from the development of model theory. 
As the field of model theory has developed over the last century, its methods have allowed for 
serious interaction with algebraic geometry, number theory, and general topology.  
One of the unique
successes of model theoretic classification theory   
is precisely in allowing for a kind of distillation and focusing of the information gleaned from the opening up of the hierarchy of infinities 
into definitions and tools for studying specific mathematical structures. 

By the time they are used, the form of the tools may not reflect this influence.  However, if we are interested in 
advancing further, 
it may be useful to remember it. 
Along this seam, so to speak, may be precisely where we will want to re-orient our approach. 

\br

\subsection*{The model theoretic point of view} 
 \blank 

\br

The model theoretic setup is designed to allow in a specific way for placing a given infinite mathematical structure in a class or family in which size 
and certain other features, such as the appearance of limit points, 
may vary.

Suppose we wanted to look abstractly at the structure the reals carry when viewed as an ordered field.  
We might consider $\mathbb{R}$ as a set decorated by the following:   
a directed graph edge representing the ordering, 
a first directed hyperedge 
representing the graph of addition, a second directed hyperedge representing the graph of multiplication, 
perhaps a third for the graph of subtraction, and two constant symbols marking the additive and multiplicative identities.  Stepping back, and retaining only this data of a set of size continuum along with the data of 
which sets of elements or tuples correspond to which constants, edges, or hyperedges, we might try 
to analyze the configurations which do or do not arise there.  

Suppose we were interested in the structure the group operation gives to the discrete Heisenberg group $\mathcal{H}$.  
We might consider $\mathcal{H}$ simply as an infinite set along with the data of the multiplication table. 
A priori, this setup just records a countably infinite set made into a group in the given way; it 
doesn't a priori record that its elements are matrices, much less uni-upper-triangular matrices over $\mathbb{Z}$.

\br

These examples suggest how models arise -- simply as sets decorated by the data of relations or functions we single out 
for study.\footnote{%Note that the edges and hyperedges mentioned above are simply two- and three-place relations. That is, a
A $k$-place relation \label{here}
on a set $X$ is a subset of $X^k$. The set of first coordinates of a binary relation is called its domain. A $k$-place function on a set $X$ is a single-valued 
binary relation whose domain is $X^k$. To specify a model, we first choose a language, which can include relation symbols, function symbols, and 
constant symbols. Then a model is given by the following data: a set $X$, called the domain;  for each relation symbol, a relation on 
$X$ of the right arity; for each function symbol, a function on $X$ of the right arity; for each constant symbol, an element of $X$.}

The initial loss of information in such a representation  
will be balanced by the 
fact that it allows us to place a model within a class %\emph{elementary class}
 and to study models in the class alongside each other. 
From the model theoretic point of view -- the following statement is a starting point for investigation, not its conclusion --  
this class contains all other models which differ from $M$ in inessential ways.

To place our model $M$ in its class, we consider the \emph{theory} of the model, that is,  
the set of all sentences of first order logic 
which hold in $M$. 
The elementary class of $M$ is the class of all other models with the same theory. 
We may frame our study as: of such classes, or of theories. 

Every model carries what are essentially derived relations, the boolean algebras of definable sets (see the Appendix). 
We might say very informally that the theory of a model takes a photograph of these boolean algebras  
which remembers only finitary information, such as which finite intersections are and are not empty. 
The models which share the theory of $M$ 
will also share this photograph. 
 
The differences among models with the same theory come essentially from the infinite intersections 
which are finitely approximated, in other words, with the filters and ultrafilters on the boolean algebra of definable sets, as we now explain.  

\br

\subsection*{Limit points}  
 \blank 

\br

To see what engenders variation within an elementary class, the following story
will motivate the definition of \emph{type}.  
(The motivation is in the telling, not an historical assertion. Filters independently arose in other, earlier 
contexts in the early part of the twentieth century.)

Writing their book on general topology in 1937, Bourbaki were discussing whether the definition of 
limit could be liberated from the countable. Cartan's suggestion was effectively that working in a topological space $X$, 
one might turn the problem around and look from the point of view of a limit point. Given a point $x \in X$, the family of neighborhoods containing $x$ 
had certain very nice features -- such as upward closure, closure under finite intersection -- which may be abstracted as 
the definition of \emph{filter}.\footnote{I learned this story from Maurice Mashaal's biography of Bourbaki, 
which 
also cites biographical work of Liliane Beaulieu. Regarding Cartan's definition of filters:  ``At first [the others] 
met the idea with skepticism, but 
Chevalley understood the importance of Cartan's suggestion and even proposed another idea 
based on it (which became the concept of ultrafilters). Once the approval was unanimous, someone yelled 
<< boum ! >> (French for ``bang!'') to announce that a breakthrough had been made -- this was one of Bourbaki's many 
customs.'' \cite{mashaal}}

\begin{defn} \label{d:filter}
For $I$ an infinite set, 
$\mcf \subseteq \mcp(I)$ is a \emph{filter} on $I$ when 
(i) $A \subseteq B \subseteq I$ and $A \in \mcf$ implies $B \in \mcf$, (ii) $A, B \in \mcf$ implies $A \cap B \in \mcf$, 
(iii) $\emptyset \notin \mcf$. 
\end{defn}

Conversely, to any filter, one can assign a (possibly empty) set of limit points: those elements of $I$ which belong to 
all $A \in \mcf$.  In defining a filter, we may restrict to any boolean algebra $\mathcal{B} \subseteq \mcp(I)$, 
asking that $\mcf \subseteq \mathcal{B}$, and adding that $B$ in item (ii) belong to $\mathcal{B}$. 
In the model theoretic context, 
this idea gives us a natural way to define limit points for any model, not requiring a metric or a topology per se: 

\begin{defn} \label{d:pt}
Informally, a \emph{partial type $p$ over a model $M$} is a filter on $M$ for the boolean algebra of $M$-definable subsets of $M$. % of elements. 
[More correctly, it is a set of formulas with parameters from $M$, whose solution sets in $M$ form such a filter.\footnote{This difference is visible in the idea of realizing a given type over $M$ in a model $N$ extending $M$. The new limit point will belong to the solution sets, 
in $N$, of the formulas in the type. Note that \ref{d:pt} describes types of \emph{elements}, corresponding to sets of formulas in one free variable plus parameters from $M$. For each $n>1$ there is an analogous $\ts_n(M)$ describing types of $n$-tuples.}]
It is a \emph{type} if it is maximal, i.e. not strictly contained in any other partial type over $M$. The Stone space $\ts(M)$ is the set of types over $M$.
\end{defn} 

A type is \emph{realized} if a corresponding limit point exists in the model, otherwise it is \emph{omitted}.  For example, if 
$M = (\mathbb{Q}, <)$ is the rationals considered as a linear order, $\ts(M)$ includes a type for each element of $M$, which are realized, 
and types for each irrational Dedekind cut, for $+\infty$, for $-\infty$,  and for various infinitessimals, which are omitted. 

The \emph{compactness theorem} for first order logic 
ensures that for any model $M$ and any type or set of types over $M$ 
we can always find an extension of $M$ to a larger model in the same class in which all these types are realized. (Put otherwise, we may realize types 
without changing the theory.)  
Types are fundamental objects in all that follows. From the depth and subtlety of their interaction comes 
much of the special character of the subject.\footnote{And a certain possible conversion of combinatorial into algebraic information, 
as in the remarkable group configuration theorems of Zilber and Hrushovski.}

In Definition \ref{d:pt}, we may use ``$A$-definable sets'' instead of ``$M$-definable sets'' for some $A \subseteq M$ 
[more correctly, formulas in one free variable with parameters from $A$]. 
In this case, call $p$ a type or partial type \emph{of $M$ over $A$}.
Then the following fundamental definition, from work of Morley and of Vaught in the early 1960s,  
generalizing ideas from Hausdorff on $\eta_\alpha$ sets,
gives a measure of the completeness of a model. 

\begin{defn}  \label{d:sat} 
For an infinite cardinal $\kappa$, we say a model $N$ is \emph{$\kappa^+$-saturated} if every type of $N$ over every 
$A \subseteq N$ 
% over every submodel $M$ of $N$ 
of size $\leq \kappa$ is realized in $N$. 
\end{defn}

For example, if $M$ is an algebraically closed field then for any subfield $K$ of $M$, the types of $M$ over $K$ will include a distinct type for each minimal 
polynomial over $K$ (describing a root) and one type describing an element transcendental over $K$. Since it is algebraically closed, 
$M$ will be $\kappa^+$-saturated if 
and only if it has transcendence degree at least $\kappa^+$.

\br

\subsection*{Towards classification theory}  
 \blank 

\br

Having described a model theoretic point of view -- first, regarding a given mathematical object as a model; 
second, placing it within an elementary class of models sharing the same first-order theory; third, studying as our 
basic objects these theories, looking both at how models may vary for a given theory (by paying close attention to the structure of 
types) and at structural differences across theories -- some first notable features of this setup are: 

\sbr
\begin{enumerate}[a)]
\item \emph{one can study the truth of statements of first order logic by `moving' statements among models in allowed ways}.\footnote{ A clever example is 
Ax's proof that any injective polynomial map from $\mathbb{C}^n$ to $\mathbb{C}^n$ is surjective. For each 
finite $k$ and $n$, there is a sentence $\vp_{n,k}$ of first order logic in the language $\{ +, \times, 0, 1 \}$ 
asserting that any injective map given by $(x_1,\dots, x_n) \mapsto (p_1(x_1,\dots, x_n), \dots, p_n(x_1,\dots, x_n))$ where the $p_i$ are of degree $\leq k$ is surjective. 
For each prime $p$, we may 
write $\bar{F}_p$ as the union of an increasing chain of finite fields. The assertion $\vp_{n,k}$ is true in each finite field 
because it is finite, and it follows (e.g. from its logical form) that $\vp_{n,k}$ holds in $\bar{F}_p$. For $\de$ any nonprincipal 
ultrafilter on the primes, $\prod_{p \in P} \bar{F}_p/\de$ is isomorphic to $\mathbb{C}$, so by \lost theorem 
 $\vp_{n,k}$ is true in $\mathbb{C}$.} 

\sbr
\item \emph{when working in a given elementary class, 
unusual constraints observed on variation of models may give leverage for a structural understanding of all models in the class}.  
E.g. Morley's theorem, the `cornerstone' of modern model theory, says that if a countable theory has only one model up to 
isomorphism in \emph{some} uncountable size $\kappa$, 
then this must be true in every uncountable size. Moreover its models must behave analogously to 
algebraically closed fields of a given characteristic in the sense that, e.g., there are prime models over sets, there are relatively few types, and for each model there is a single invariant, a dimension (the equivalent of transcendence degree) giving the isomorphism type. 

\sbr
\item[c)] \emph{simply understanding the structure of the definable sets, say in specific classes containing 
examples of interest, can already involve deep mathematics}. 
For instance, Tarski's proof of quantifier elimination for the reals and the cell decomposition theorem for o-minimal structures.

\end{enumerate}
In the examples given so far, as is often the case in mathematics, the specific role of the infinite may be mainly as a kind of foil reflecting 
the fine structure of compactness, irrespective of the otherwise depth of proofs.

To see the interaction of model theory and set theory which we invoked at the beginning, 
we need to go further up and further in. 
(As an aside, already in 
Hilbert's remarks, via Church, there is an implicit parallel between the   
understanding of infinite sizes and the development 
of different models.\footnote{  \hspace{2mm}``Hilbert  does not  say  that  the  
order   in  which  the [list of 23] problems  are  numbered  gauges  their    relative    
importance,  and  it  is  not  meant  to  suggest  that  he  intended  this.  
But  he  does  mention  the  arithmetical  formulation   of  the  concept  
of  the  continuum  and  the  discovery   of  non-Euclidean  geometry  as  
being  the  outstanding  mathematical  achievements   of  the  preceding  
century,  and  gives  this  as  a  reason  for  putting  problems  in  these  
areas first, '' \cite{church}.})  
Suppose we step back and study the class of  all theories. 

A thesis of Shelah's groundbreaking \emph{Classification Theory} (1978) is that one can find 
dividing lines among the class of first-order theories.  A dividing line 
marks a sea change in the 
combinatorial 
structure. 
(The assertion that something is a dividing line requires evidence on both sides:
showing that models of theories on one side are all complex in some given sense, while models of theories on the other 
side admit some kind of structure theory.) 
A priori it is not at all obvious that these should exist.  
Why wouldn't the seemingly unconstrained range and complexity of theories allow for some kind of continuous gradation 
along any reasonable axis? The example of extremal combinatorics may give a hint: graphs of a given large 
finite size 
 are not so easily classifiable, but by examining asymptotic growth rates of certain phenomena, jumps may appear. 

Stability, the dividing line which has most profoundly influenced the present field, arises in \cite{Sh:a} from counting limit points. 
For a theory $T$ and an 
infinite cardinal $\lambda$, we say $T$ is \emph{$\lambda$-stable} if for every model $M$ of $T$ of size $\lambda$, $|\ts(M)| = \lambda$. 
Conversely, if some model $M$ of $T$ of size $\lambda$ has $|\ts(M)| > \lambda$, $T$ is \emph{$\lambda$-unstable}. For a given $\lambda$, 
$T$ is either $\lambda$-stable or $\lambda$-unstable by definition. But varying $\lambda$, the gap appears\footnote{$|T|$, the size of $T$, will be the maximum of the size of the language and $\aleph_0$.}:

\begin{theorem}[Shelah 1978]  \label{th:stable}
Any theory $T$ is either \emph{stable}, meaning stable in all $\lambda$ such 
that $\lambda^{|T|} = \lambda$, or \emph{unstable}, 
meaning unstable in all $\lambda$.
\end{theorem}

This theorem materializes in step with the development of the internal 
structure theory. 
The set theoretic scaffolding is not only in the statement, but intricately connected to its development.    
A few examples from \cite[II-III]{Sh:a} will give a flavor: 
\begin{enumerate}[i)]

\item it turns out stability is local: if $T$ is unstable,  
then there is a single formula $\vp$ such that 
in all $\lambda$, we can already get many types just using definable sets which are instances of $\vp$. (This leads to discovering instability has a characteristic combinatorial configuration, the order property.)
Its proof is  
a counting argument relying on the fact that $T$ is unstable in some $\lambda$ such that $\lambda = \lambda^{|T|}$. 

\item a characteristic property of stable theories is that once there is enough information, types have unique generic extensions to types over larger sets. 
This is first explained by the finite equivalence relation theorem  
which studies types over sets $A$ of size at least $2$ in models which are $(|A|^{|T|})^+$-saturated.

\item conversely, large types are essentially controlled by 
their restrictions to ``small'' sets. The cardinal defining this use of ``small,''  $\kappa(T)$, is $\leq |T|^+$.
Above this cardinal the mist clears and stability's effects are easier to see; one can e.g. 
characterize larger saturated models of stable theories as those which are $\kappa(T)$-saturated and 
every maximal indiscernible set has the cardinality of the model. 

\item instability has a more random form and a more rigid form, and at least one must occur.  To see the difference between the two by counting types, one has to know that $2^\lambda$ 
can potentially be different from the number of cuts in a dense linear order of size $\lambda$ (or the number of branches in a tree with $\lambda$ nodes), which was originally noted by appeal to an independence result. 
\end{enumerate}

\br

Of course, to say infinite cardinals are strongly connected to the development of stability doesn't mean they are necessarily there at the end. 
The order property, definability of types, forking, the independence property, and the strict order property, to name a few, 
don't bear the imprint of their origin. This translation is part of model theory's power.

Looking forwards: the effect of stability, inside and outside model theory, has been significant.  Despite this conclusive evidence that 
some dividing lines do exist, and that they can be very useful, further ones have been challenging to find. We know very few in the vast territory of unstable theories, found -- like stability -- one by one in response to specific counting problems. 
To go further, perhaps we can try to shift the way in which 
set theory sounds out model theoretic information.

The reader may wonder: is model theory being described as a kind of extension to the infinite setting of extremal arguments in combinatorics, with the hierarchy of 
infinite cardinals replacing the natural numbers?  
This analogy is challenging, but incomplete. It is incomplete 
because the finitary, extremal picture doesn't seem to provide a precedent or explanation for the role 
of model theory, 
which builds in a remarkable way a bridge between the infinite combinatorial world and a more algebraic one. Still, it is 
challenging because it leads us to ask what in the infinite setting may play the role of 
those crucial tools of the combinatorial setting, which may seem to have little 
place in current model theoretic arguments --  namely, probability and randomness.

\br

\subsection*{Ultrapowers} \blank 

\br

Only in the move to ultrapowers 
does one really recover, albeit in a metaphorical way, 
that other key ingredient of extremal arguments, the understanding of probability and average behavior.  

Stability arises from counting limit points. Recall from \ref{d:sat} that saturation is a notion of completeness for a model: $\lambda^+$-saturation 
means all types over 
all submodels of size at most $\lambda$ are realized. 
The ultrapower construction, given formally below, starts with a given model and amplifies it -- staying within the 
elementary class -- according to a specific kind of  
averaging mechanism, an ultrafilter.  The resulting larger model, which depends only on the model we began with and the ultrafilter, 
is called an ultrapower.  The level of saturation in the ultrapower reflects whether the given averaging mechanism, 
applied to the given model, leads to the appearance of many limit points of smaller sets.  

If so, this may indicate either that the ultrafilter is powerful, or that the types of the model are not complex.  
Since we can apply the same ultrafilter to different models and 
compare the results, however, we can use this construction to compare the `complexity' of different models (and, dually, of ultrafilters). 
Restricting to the powerful class of \emph{regular} ultrafilters, whose ability to produce saturation in ultrapowers will be an invariant of the elementary class of the model we begin with, we can use this construction to compare the complexity of theories.  Informally 
for now, the relation on complete countable theories setting $T_1 \trianglelefteq T_2$
\begin{quotation}
\noindent  if for any regular ultrafilter $\de$, if $\de$-ultrapowers of models of $T_2$ are sufficiently saturated, 
so are $\de$-ultrapowers of models of $T_1$
\end{quotation} is Keisler's order, defined in 1967. It is a pre-order on theories, considered as a partial order on the equivalence classes of theories. 

The theorem which convinced the author of this essay, reading around 2005, that it was urgent to study $\trianglelefteq$ further was a
theorem in Shelah's \emph{Classification Theory} in a chapter devoted to the ordering. The theorem says that the union of the 
first two equivalence classes in Keisler's order is precisely the stable theories. 

This theorem can be understood as saying that the class of stable theories, which we can see by counting types, can also be seen by asking about 
good average behavior. 
Beyond stability, our counting is less useful, and yet the other half, about average behavior, 
retains its power.

 \br
\newpage

\subsection*{In more detail} \blank 

\br

The idea of a filter was used in \ref{d:pt} above to find limit points, but it can also be used to give averages. Maximal filters, called ultrafilters, 
can be thought of as a coherent choice of which subsets of a given set $I$ are ``large.'' 

\begin{defn}
For an infinite set $I$, $\de \subseteq \mcp(I)$ is an ultrafilter if it is a filter not strictly contained in any other filter. 
$($We will assume ultrafilters contain all co-finite sets.$)$
\end{defn}

The ultraproduct by $\de$ of a family of models $\langle M_i : i \in I \rangle$ is a model, built in two steps, reflecting the definition of model. 
First, we define the domain. Identify two elements $\langle a[i] : i \in I \rangle$, $\langle b[i] : i \in I \rangle$ of the Cartesian 
product $\prod_i M_i$ if $\{ i \in I : a[i] = b[i] \} \in \de$. Definition \ref{d:filter} makes this an equivalence relation, and the 
domain of our ultraproduct $N$ is the set of equivalence classes. Next, fix for transparency a representative of each equivalence class, so that 
for $a \in N$ and $i \in I$, ``$a[i]$'' makes sense. The relations, functions, and constants of our language are defined on the 
ultraproduct by consulting the average of the models: e.g. we say a given $k$-place relation $R$ holds on 
$a_1,\dots, a_k$ in $N$ iff $\{ i : R(a_1[i],\dots, a_k[i]) $ holds in $M_i \} \in \de$, and for an $n$-place function $f$, 
define $f(a_1,\dots, a_n)$ 
to be the equivalence class of $\langle b[i] : i \in I \rangle$ where $b[i] = f(a_1[i], \dots, a_n[i])$ computed in $M_i$. 
The special case 
of an ultrapower, where all the factor models are isomorphic,  transforms a given structure into a larger, `amplified' model in the same elementary class. 
For example, letting $\de$ be a regular ultrafilter on the set of primes, $\prod_p \overline{\mathbb{F}}_p/\de \cong \mathbb{C}$, 
but if we consider the ultrapower $M^P/\de$ where $M$ is the algebraic closure of the rationals, we also get $\mathbb{C}$.

There is a veiled interaction between the two model theoretic uses of filters: the realization of \emph{types} 
in the ultrapower, and the \emph{ultrafilter} used in the construction. This is most useful when the 
ultrafilter is \emph{regular}, ensuring that saturation depends only on finitary input from each factor 
model. For each regular $\de$, 
whether or not an ultrapower $M^I/\de$ is $|I|^+$-saturated is an invariant of the elementary class 
of the model  $M$ we began with.\footnote{That can be taken as a definition of regular; alternatively, $\de$ is a regular ultrafilter on 
$I$ if there is a set $\{ X_\alpha : \alpha < |I| \} \subseteq \de$ such that the intersection of any infinitely many of its elements is empty. Regular 
ultrafilters are easy to build and exist on any infinite set.}

Keisler's suggestion was that this could be used to compare theories. 

\begin{defn}[Keisler's order, 1967] Let $T_1, T_2$ be complete countable theories. 
\[ T_1 \trianglelefteq T_2 \]
if for every infinite $\lambda$, every regular ultrafilter $\de$ on $\lambda$, every model $M_1$ of $T_1$ and every model 
$M_2$ of $T_2$, if $(M_2)^\lambda/\de$ is $\lambda^+$-saturated, then $(M_1)^\lambda/\de$ is $\lambda^+$-saturated. 
\end{defn}

Informally, say ``$\de$ saturates $T$'' if for some (by regularity of $\de$, the choice 
does not matter) model $M$ of $T$, $M^\lambda/\de$ is $\lambda^+$-saturated: all limit points over small submodels appear. 
Then Keisler's order puts $T_1$ less than $T_2$ if 
every regular ultrafilter able to saturate $T_2$ is able to saturate $T_1$.  Note that any two theories may in principle be compared -- 
algebraically closed fields of fixed characteristic, random graphs, real closed fields. 
Keisler proved his order was well defined and had a minimum and a maximum class (he gave a sufficient condition for membership in each), and asked about its 
structure.

The crucial chapter on Keisler's order in \cite{Sh:a} was already mentioned. Its structure on the unstable, non-maximal theories was 
left there as an important open question. 
Following \cite{Sh:a}, work on Keisler's order stalled for about thirty years. The question was reopened in Malliaris' thesis and in the 
series of papers \cite{mm1}, \cite{mm2}, \cite{mm3}, \cite{mm4}, \cite{mm5}, guided by the perspective described above.  
Then in the last few years, a very productive ongoing collaboration of 
Malliaris and Shelah  
\cite{MiSh:996}, \cite{MiSh:997}, \cite{MiSh:999}, \cite{MiSh:E74}, \cite{MiSh:998}, \cite{MiSh:1030}, %\cite{MiSh:1051}, 
\cite{MiSh:1050}... has advanced things considerably.

\br

\subsection*{The current picture} \blank 

\br

Along this road, what does one find? 

\sbr
\noindent (a) First, we do indeed see evidence of dividing lines --  many more than previously thought.  Much remains to be done in understanding them 
and in characterizing the structure/nonstructure which come with the assertion of a dividing line, but already their appearance, in 
a region of theories thought to be relatively tame, is surprising and exciting. 
In \cite{MiSh:1050}, Malliaris and Shelah prove that Keisler's order has infinitely many classes. The theories which 
witness these different classes come from higher analogues of the countable triangle-free random graph, originally studied by Hrushovski \cite{hrushovski1}: 
the infinite generic tetrahedron-free three-hypergraph, the infinite generic 4-uniform hypergraph with no complete hypergraph on 5 vertices, 
and so on.  The proof shows they 
may have very different average behavior, 
as reflected in their differing sensitivity to a certain degree of calibration in the ultrafilters.  

These results build on advances in ultrafilter construction, which allow for a greater use of properties of 
cardinals, even for ultrafilters in ZFC. 

Several incomparable classes are known \cite{MiSh:1124}, \cite{ulrich}, and it may be that future work will reveal many.  
Perhaps the way that such averages could be perturbed or distorted, and by extension the structure of dividing lines among 
unstable theories, will be much finer than what we now see.  
If so, even independence results could be quite useful model theoretically. 
These may simply witness that the boolean algebras associated to different theories 
are essentially different, because they react differently to 
certain exotic averaging mechanisms, when these appear.  The internal theories of each equivalence class, giving an account of 
what {allows} for the different reactions, would presumably be, like IP or SOP, absolute. 

\sbr

\noindent (b) Second, this line of work has led to some surprising theorems about the finite world.  These theorems have the following general form. 
We know that among theories with infinite models, stability is a dividing line, with models of stable theories admitting a strong structure theory. 
There is a specific combinatorial configuration, the order property, which (in infinite models) characterizes instability.   
In an infinite graph, instability for the edge relation would correspond to having arbitrarily large \emph{half-graphs}, that is, for all $k$ having 
vertices $a_1,\dots, a_k$ and $b_1,\dots, b_k$ with an edge between $a_i$ and $b_j$ iff $i<j$. 
(Note that there are no assertions made about edges among the $a$'s and among the $b$'s, so in forbidding $k$-half graphs, 
we forbid a family of configurations.) 
The thesis of Malliaris and Shelah \cite{MiSh:978} is essentially that 
finite graphs with no long half-graphs, called \emph{stable graphs}, behave much better than 
all finite graphs, in the sense predicted by the infinite case. 

Szemer\'edi's celebrated regularity lemma says, 
roughly, that for every $\epsilon > 0$ there is $N = N(\epsilon)$ such that any sufficiently large graph 
may be equitably partitioned into $k \leq N$ pieces such that all but at most $\epsilon k^2$ pairs of pieces have the 
edges between them quite evenly distributed (i.e. are $\epsilon$-regular).     
The elegant picture of this lemma absorbs the general complexity of graphs in two ways: first, 
by work of Gowers, $N$ is a very large function of $\epsilon$; second, as noticed independently by 
a number of researchers, the condition that some pairs of pieces be irregular cannnot be removed, as 
shown by the example of half-graphs \cite{k-s}. As a graph theorist, one might expect half-graphs to be 
just one example of bad behavior, not necessarily unique;  but in light of the above, 
a model theorist may guess that in the absence of 
long half-graphs one will find structured behavior. The stable regularity lemma of \cite{MiSh:978} shows that indeed, half-graphs 
are the only reason for irregular pairs: finite stable graphs admit regular partitions with no irregular pairs and 
the number of pieces singly-exponential in $\epsilon$.  

A second theorem in that paper, a stable Ramsey theorem, proves that for each $k$ there is $c=c(k)$ such that 
if $G$ is a finite graph with no $k$-half graphs then $G$ contains a clique or an independent set of size 
$|G|^c$, much larger than what is predicted by Ramsey's theorem. This meets the prediction of the Erd\H{o}s-Hajnal 
conjecture, which says that for any finite graph $H$ there is $c = c(H)$ such that if $G$ contains no induced copies of $H$ 
then $G$ contains a clique or an independent set of size $|G|^c$. But very few other cases of this conjecture are known.  
What is the contribution of the infinite here? 
The infinite version of Ramsey's theorem says that a countably infinite graph contains a countable clique or 
a countable independent set. Extending this to larger cardinals doesn't get far:   
Erd\H{o}s-Rado shows the graph and the homogeneous set may not, in general, increase at the same rate. 
Model theory, however, refracts this result across different classes of theories, and across dividing lines, and in some classes, such as the stable theories, 
it behaves differently: the existence of large sets of indiscernibles in stable theories implies, 
a fortiori, that an infinite stable graph of size $\kappa^+$ will have a clique or independent set of size $\kappa^+$, 
much larger than predicted by Erd\H{o}s-Rado. Once one knows where to look, one can find the analogous phenomenon in the 
finite case (also for hypergraphs). 

The stable Ramsey theorem was applied by Malliaris and Terry \cite{malliaris-terry} 
to re-prove a theorem in the combinatorics literature, 
by re-organizing the proof into cases which take advantage of the stable Ramsey theorem, and thus to obtain 
better bounds for the original theorem; finitary model-theoretic analysis may be useful even where model theoretic 
hypotheses are not used in the theorems. 

It may seem that these theorems of \cite{MiSh:978}, from the first joint paper of Malliaris and Shelah, could in principle have been discovered earlier, 
and yet they were not.   They belong to the perspective of this program in a deep way. 
They were 
motivated by work in \cite{mm2} and \cite{mm3}, directed towards 
Keisler's order, which first 
applied Szemer\'edi regularity to study the complexity of formulas (and showed that the simple theories, of special interest in Keisler's order, 
were in some sense controlled by stable graphs). 

\sbr

\noindent (c) 
Third, by means of these methods model theory has paid an old debt to set theory and general topology, by solving a seventy-year-old problem about  
cardinal invariants of the continuum. Two infinite cardinals, $\xp$ and $\xt$, known to be uncountable but no larger than the continuum, are 
shown in Malliaris and Shelah \cite{MiSh:E74}, \cite{MiSh:998} to be unconditionally equal. The proof is model-theoretic, and 
comes in the context of the solution 
of an a priori unrelated problem, determining a new sufficient condition for maximality in Keisler's order. 

In slightly more detail, in order for a regular ultrafilter to handle the most complex theories, those in the maximum class, it must be in some sense very balanced. Distortions and 
so to speak imperfections which might pass unnoticed in more robust   
theories will translate immediately in maximal theories to the omission of types. 
However, a surprising fact from \cite{Sh:a} is that what is needed for the theory to be complex is not necessarily that it be expressive. 
A kind of brittleness or overall rigidity as exemplified by the theory of linear order will also suffice. 
Remarkably, it turns out even less will suffice: the engine of the proof in \cite{MiSh:998} is in showing that
if the ultrafilter can ensure certain paths through trees have upper bounds, it must be strong enough to produce the needed 
limit points for any theory. 
The resulting comparison of theories whose models contain the relevant 
trees, to models of linear orders, via ultrapowers, turns out to be parallel in a precise sense to the comparison of $\xp$ and $\xt$. 
It was possible to give a fundamental model-theoretic framework encompassing both problems and so to solve them both. 
A commentary and an expository account of the proof are \cite{moore}, \cite{malliaris-casey}.   

Still, a model theoretic necessary condition for maximality remains open. 

\br

\noindent
It has been almost ten years since \cite{mm1}.  
Profound questions remain, urgent in their simplicity.

\newpage
\subsection*{Appendix: on definable sets}

\blank

\br

By convention, we always assume our language $\ml$ contains a binary relation symbol $=$, 
and that in every $\ml$-model this symbol is interpreted as equality. 
Besides the symbols of $\ml$, our alphabet for building formulas includes 
infinitely many variables along with logical symbols $(, ), \land, \lor, \neg, \iff, \implies, \forall, \exists$.

For awhile let $M_\star$ denote a model for the language $\ml_\star = \{ +, \times, -,  0, 1, < \}$, where $+$, $\times$, $-$ 
are binary function symbols, 
$<$ is a binary relation, and $0, 1$ are constants.    
Let us say the domain of $M_\star$ is $\mathbb{R}$ and the symbols have their usual interpretation.\footnote{Pedantically, the form of the symbols makes no demands on their interpretation (other than the basic conditions in the footnote on p. \pageref{here}); we could, for example, give a perfectly valid model by interpreting both $+$ and $\times$ by projection onto the first coordinate.} 

The \emph{terms} of a language are elaborate names. 
We define terms by induction. 
All variables and constant symbols are terms; 
 if $f$ is a $k$-place function symbol and $t_1,\dots, t_k$ are terms, then $f(t_1,\dots, t_k)$ is a term; and a finite string of allowed symbols 
is a term iff it can be built in finitely many steps in this way. 

The $\ml_\star$-terms include $0$, $1$, $1+1$ which we may abbreviate $2$, $x+1$, 
$x \times y$, $(x \times 1) + (y + 0)$, $((x \times x) \times x)$ which we may abbreviate $x^3$. Using similar abbreviations, and 
dropping parentheses for readability, $x^5 + 15 x^2 + 3x + 5$ is also a term. 

A key feature of terms is that 
\emph{if} we are working in a model, and we are given a term along with instructions of which elements of the model 
to put in for which, if any,  
variables in the term (recall that in any model, any constant symbols must already refer to specific elements), then 
the term will evaluate unambiguously to some other element of the model. 

Next, by induction, we define \emph{formulas}. Atomic formulas are assertions that a given relation symbol of our language holds 
on a given sequence of terms.  (In $\ml_\star$, the relation symbols are $=$ and $<$, so 
these will include  $x^3 + 5x + 2 = 0$ and also $5 + x + 15y^2 > 37 - z$.)
Atomic formulas are formulas. 
If $\vp$ is a formula, then $\neg(\vp)$ is a formula. If $\vp$ and $\psi$ are formulas, then $(\vp \land \psi)$, 
$(\vp \lor \psi)$, $(\vp \implies \psi)$, $(\vp \iff \psi)$ are also formulas. If $\vp$ is a formula, and $x$ is a variable, then 
$(\exists x) \vp$ and $(\forall x) \vp$ are formulas. A finite string of allowed symbols is a formula iff it can be built in finitely many 
steps in this way.  

In our example here are some more formulas:   $(\forall x)(x+0 = x)$, $(\exists x) ( y + x^2  = z)$. 
Note an important difference between the two. $(\forall x)(x+0 = x)$ is an assertion which will be true or false in any given model; 
in our given $M_\star$ it is true. By contrast, $(\exists x) ( y + x^2  = z)$ is neither true nor false, since it has 
two free variables; rather, it has a \emph{solution set}, 
the pairs $(a,b)$ of elements of $M$ such that $(\exists x)(a + x^2 = b)$.  In any model $N$, the 
solution sets of formulas with one or more [but always finitely many] 
free variables are called the \emph{definable sets}. 
The closure properties of the set of formulas show that for each $n$, the definable sets on $N^n$ form a boolean algebra. 

The formulas with no free variables are called \emph{sentences}, and the theory of a model $N$ is the set of all sentences 
which hold in $N$. The elementary class of $N$ is the class of all other models in the same language with the same theory.
The reason a theory may make assertions about definable sets which are meaningful across different models is by 
referring to their defining formulas.  

%We need one more definition. 
When a formula has many free variables, it may be useful to look at the restricted solution set 
we get after specifying that certain of the free variables take certain values in the model. For example, 
in the formula $x y^2 + z y + w = 0$ with free variables $x,y,z,w$, we might want to consider the solution set under specific values of $x, z$ and $w$. 
Such a solution set is called \emph{definable with parameters}, the specific values being the parameters.  
We may wish to record their provenance: given a subset $A$ of a model $N$, 
the sets definable with parameters from $A$ are called \emph{$A$-definable sets}. 
Finally, a word on types. Given $M$ and $A \subseteq M$, 
%the $A$-definable subsets of $M$ form a boolean algebra. 
the set of formulas with parameters from $A$ and (say) one free variable can be made into a boolean algebra 
once the formulas are identified up to logical equivalence (equivalently, identified if 
they define the same set in $M$). Its Stone space is the set of types of $M$ over $A$ in the sense of 
\ref{d:pt} above, and its compactness as a topological space is explained by the compactness theorem.  

Some examples -- here, definable means with parameters: 

\begin{enumerate}[1)]
\item  in the model $M_\star$ above, the definable sets include the semialgebraic sets (and it is a theorem that they are exactly the 
semialgebraic sets).  Its elementary class is the class of real closed ordered fields. 

\item if $\ml = \{ +, \times, -, 0, 1 \}$ and $M$ is the algebraic closure of the rationals on which the symbols have their usual interpretation, 
the definable sets include (and, in fact, are) the constructible sets. 
The elementary class of $M$ is the class of all algebraically closed fields of characteristic zero. 

\item if $\ml = \{ < \}$, and $M$ is the rationals on which $<$ has its usual interpretation, the definable sets in one free variable 
are finite unions of points and intervals (and the definable sets in $k >1$ free variables satisfy a 
cell decomposition theorem).  The elementary class of $M$, the class of  
dense linear orders without a first or last element, has only 
one countable model, up to isomorphism, but $2^\lambda$ nonisomorphic models of each uncountable size $\lambda$. 
\end{enumerate}

\vspace{8mm}

\noindent \emph{Acknowledgments}.  B. Mazur, A. Peretz, and S. Shelah made helpful comments.  
The Institute for Advanced Study contributed some wonderful working conditions. 
The National Science Foundation provided partial support through grant 1553653 and grant 
1128155 to the IAS. Thank you.

\br

\end{document}